     \documentstyle[12pt,fleqn]{article}
\expandafter\chardef\csname pre amssym.def at\endcsname=\the\catcode`\@
\catcode`\@=11
\def\hexnumber@#1{\ifcase#1 0\or 1\or 2\or 3\or 4\or 5\or 6\or 7\or 8\or
 9\or A\or B\or C\or D\or E\or F\fi}
\font\tenmsa=msam10
\font\sevenmsa=msam7
\font\fivemsa=msam5
\newfam\msafam
\textfont\msafam=\tenmsa
\scriptfont\msafam=\sevenmsa
\scriptscriptfont\msafam=\fivemsa
\edef\msafam@{\hexnumber@\msafam}
\def\emptybox{\mathrel{\mathchar"0\msafam@03}}
\font\tenmsb=msbm10
\font\sevenmsb=msbm7
\font\fivemsb=msbm5
\newfam\msbfam
\textfont\msbfam=\tenmsb
\scriptfont\msbfam=\sevenmsb
\scriptscriptfont\msbfam=\fivemsb
\edef\msbfam@{\hexnumber@\msbfam}
\def\ltimes{\mathrel{\mathchar"0\msbfam@6E}}
\def\rtimes{\mathrel{\mathchar"0\msbfam@6F}}
\catcode`\@=\csname pre amssym.def at\endcsname
\newtheorem{theorem}{\bf Theorem}[subsection]
\newtheorem{proposition}[theorem]{\bf Propositon}

\newtheorem{lemma}[theorem]{\bf Lemma}

\newtheorem{example}{\bf Example}[subsection]
\newenvironment{proof}{\par\noindent{\bf Proof.}}%
{\kern 20pt $\emptybox$ \medskip}

\def\subsec{\subsection{}\vskip -22pt\hskip \parindent}

\def\cO#1{{{\cal O}(#1)}}
\def\RR{I\!\!R}

\def\ZZ{Z\!\!\!Z}
\def\NN{I\!\!N}
\def\T#1{{\rm T}(#1)}
\def\V{{\rm V}}
\def\C{{\cal C}}
\def\op{{\rm op}}
\def\id{{\rm id}}
\def\Obj{{\rm Obj}}

\def\krr{\kern -.16667em}%
\def\kr{}%
\def\krrr{\kern -.3\unitlength}%
\def\krl{}%
\def\krrrr%
{\kern -5\unitlength}
\newlength{\textwd}%
\def\hhstep{\kr\kr
\kern -.5\unitlength}
\def\hstep{\kr\kr
\kern .5\unitlength}
\def\step{\kr\kr
\kern \unitlength}
\def\Step{\kr\kr
\kern 2\unitlength}
\def\vvbox#1{{\offinterlineskip\vcenter{%
\def\coev{\kr
\begin{picture}(2,2)\put(1,0){\oval(2,2)[t]}\end{picture}}
\def\ev{\kr
\begin{picture}(2,2)\put(1,2){\oval(2,2)[b]}\end{picture}}
\def\hcoev{\kr
\begin{picture}(1,2)\put(.5,0){\oval(1,1)[t]}\end{picture}}
\def\hev{\kr
\begin{picture}(1,2)\put(.5,2){\oval(1,1)[b]}\end{picture}}
\def\COEV{\kr
\begin{picture}(2,2)\put(3,0){\oval(6,6)[t]}\end{picture}}
\def\EV{\kr
\begin{picture}(2,2)\put(3,2){\oval(6,6)[b]}\end{picture}}
\def\unit{\kr
\begin{picture}(0,2)
\put(0,0){\line(0,1){1}}\put(0,1.2){\circle{0.4}}
\end{picture}}
\def\counit{\kr
\begin{picture}(0,2)
\put(0,1){\line(0,1){1}}\put(0,.8){\circle{0.4}}
\end{picture}}
\def\Q##1{\kr
\begin{picture}(0,2)
\put(0,0){\line(0,1){0.4}}\put(0,1){\circle{1.2}}
\put(-0.6,0.4){\makebox(1.2,1.2)[cc]{$\scriptstyle ##1$}}
\end{picture}}
\def\O##1{\kr
\begin{picture}(0,2)
\put(0,0){\line(0,1){0.4}}\put(0,1.6){\line(0,1){0.4}}\put(0,1){\circle{1.2}}
\put(-0.6,0.4){\makebox(1.2,1.2)[cc]{$\scriptstyle ##1$}}
\end{picture}}
\def\S{\O{S}}                 \def\SS{\O{S^-}}
\def\tS{\O{\overline S}}     \def\tSS{\O{\overline S^-}}
\let\P\O
\def\dash##1{\kr
\begin{picture}(2,2)
\put(-.5,0){\dashbox{.1}(3,2){$\scriptstyle ##1$}}
\end{picture}}
\def\Dash##1{\kr
\begin{picture}(2,2)
\put(-1,0){\dashbox{.1}(4,2){$\scriptstyle ##1$}}
\end{picture}}
\def\x{\kr
\begin{picture}(2,2)
\put(0,2){\line(1,-1){2}}\put(0,0){\line(1,1){.7}}\put(2,2){\line(-1,-1){.7}}
\end{picture}}
\def\xx{\kr
\begin{picture}(2,2)
\put(0,2){\line(1,-1){.7}}\put(0,0){\line(1,1){2}}\put(2,0){\line(-1,1){.7}}
\end{picture}}
\def\hx{\kr
\begin{picture}(1,2)
\put(0,2){\line(1,-2){1}}\put(0,0){\line(1,2){.35}}\put(1,2){\line(-1,-2){.35}}
\end{picture}}
\def\hxx{\kr
\begin{picture}(1,2)
\put(0,2){\line(1,-2){.35}}\put(0,0){\line(1,2){1}}\put(1,0){\line(-1,2){.35}}
\end{picture}}
\def\d{\kr
\begin{picture}(1,2)\put(0,2){\line(1,-2){1}}\end{picture}}
\def\dd{\kr
\begin{picture}(1,2)\put(0,0){\line(1,2){1}}\end{picture}}
\def\hd{\kr
\begin{picture}(1,2)
\put(0,2){\line(1,-2){.5}}
\put(.5,1){\line(0,-1){1}}
\end{picture}}
\def\hdd{\kr
\begin{picture}(1,2)
\put(1,2){\line(-1,-2){.5}}
\put(0,1){\line(0,-1){1}}
\end{picture}}
\def\ld{\kr
\begin{picture}(1,2)
\put(1,0){\oval(2,2)[lt]}\put(1,0){\line(0,1)2}
\end{picture}}
\def\Ld{\kr
\begin{picture}(2,2)
\put(2,0){\oval(4,2)[lt]}\put(2,0){\line(0,1)2}
\end{picture}}
\def\cd{\kr
\begin{picture}(2,2)
\put(1,0){\oval(2,2)[ct]}\put(1,1){\line(0,1)1}
\end{picture}}
\def\hdcd{\kr
\begin{picture}(1,2)
\put(0,2){\line(1,-2){.5}}
\put(.5,0){\oval(1,1)[ct]}\put(.5,.5){\line(0,1){.5}}
\end{picture}}
\def\hddcd{\kr
\begin{picture}(1,2)
\put(1,2){\line(-1,-2){.5}}
\put(.5,0){\oval(1,1)[ct]}\put(.5,.5){\line(0,1){.5}}
\end{picture}}
\def\hcd{\kr
\begin{picture}(1,2)
\put(.5,0){\oval(1,1)[ct]}\put(.5,.5){\line(0,1){1.5}}
\end{picture}}
\def\Cd{\kr
\begin{picture}(4,2)
\put(2,0){\oval(4,2)[ct]}\put(2,1){\line(0,1)1}
\end{picture}}
\def\rd{\kr
\begin{picture}(1,2)
\put(0,0){\oval(2,2)[rt]}\put(0,0){\line(0,1)2}
\end{picture}}
\def\Rd{\kr
\begin{picture}(2,2)
\put(0,0){\oval(4,2)[rt]}\put(0,0){\line(0,1)2}
\end{picture}}
\def\lu{\kr
\begin{picture}(1,2)
\put(1,2){\oval(2,2)[lb]}\put(1,0){\line(0,1)2}
\end{picture}}
\def\Lu{\kr
\begin{picture}(2,2)
\put(2,2){\oval(4,2)[lb]}\put(2,0){\line(0,1)2}
\end{picture}}
\def\cu{\kr
\begin{picture}(2,2)
\put(1,2){\oval(2,2)[cb]}\put(1,0){\line(0,1)1}
\end{picture}}
\def\hdcu{\kr
\begin{picture}(1,2)
\put(1,0){\line(-1,2){.5}}
\put(.5,2){\oval(1,1)[cb]}\put(.5,1){\line(0,1){.5}}
\end{picture}}
\def\hddcu{\kr
\begin{picture}(1,2)
\put(0,0){\line(1,2){.5}}
\put(.5,2){\oval(1,1)[cb]}\put(.5,1){\line(0,1){.5}}
\end{picture}}
\def\hcu{\kr
\begin{picture}(1,2)
\put(.5,2){\oval(1,1)[cb]}\put(.5,0){\line(0,1){1.5}}
\end{picture}}
\def\Cu{\kr
\begin{picture}(4,2)
\put(2,2){\oval(4,2)[cb]}\put(1,0){\line(0,1)1}
\end{picture}}
\def\ru{\kr
\begin{picture}(1,2)
\put(0,2){\oval(2,2)[rb]}\put(0,0){\line(0,1)2}
\end{picture}}
\def\Ru{\kr
\begin{picture}(2,2)
\put(0,2){\oval(4,2)[rb]}\put(0,0){\line(0,1)2}
\end{picture}}
\def\k{\kr
\begin{picture}(1,2)
\put(0,2){\oval(2,1)[rb]}
\put(0,0){\oval(2,1)[rt]}
\put(0,0){\line(0,1)2}
\end{picture}}
\def\kk{\kr
\begin{picture}(1,2)
\put(1,2){\oval(2,1)[lb]}
\put(1,0){\oval(2,1)[lt]}
\put(1,0){\line(0,1)2}
\end{picture}}
\def\ro##1{\kr
\begin{picture}(2,2)
\put(.4,0){\oval(.8,.8)[lt]}\put(1.6,0){\oval(.8,.8)[rt]}
\put(1,0.4){\circle{1.2}}
\put(0.4,-0.2){\makebox(1.2,1.2)[cc]{$\scriptstyle ##1$}}%
\end{picture}}
\def\coro##1{\kr
\begin{picture}(2,2)
\put(.4,2){\oval(.8,.8)[lb]}\put(1.6,2){\oval(.8,.8)[rb]}
\put(1,1.6){\circle{1.2}}
\put(0.4,1){\makebox(1.2,1.2)[cc]{$\scriptstyle ##1$}}%
\end{picture}}
\def\Ro##1{\kr
\begin{picture}(4,2)
\put(1.4,0){\oval(2.8,1.2)[lt]}\put(2.6,0){\oval(2.8,1.2)[rt]}
\put(2,.6){\circle{1.2}}
\put(1.4,0){\makebox(1.2,1.2)[cc]{$\scriptstyle ##1$}}%
\end{picture}}
\def\coRo##1{\kr
\begin{picture}(4,2)
\put(1.4,2){\oval(2.8,1.2)[lb]}\put(2.6,2){\oval(2.8,1.2)[rb]}
\put(2,1.4){\circle{1.2}}
\put(1.4,.8){\makebox(1.2,1.2)[cc]{$\scriptstyle ##1$}}%
\end{picture}}
\def\r{\ro{\cal R}}              \def\rr{\ro{{\cal R}^-}}
            \def\rrr{\ro{{\cal R}^{\tilde{}}}}
\def\ra{\ro{{\cal R}_A}}        \def\rra{\ro{{\cal R}^-_A}}
\def\rb{\ro{{\cal R}_B}}        \def\rrb{\ro{{\cal R}^-_B}}
\def\rh{\ro{{\cal R}_H}}
\def\R{\Ro{\cal R}}           \def\RR{\Ro{{\cal R}^-}}
\def\Ra{\Ro{{\cal R}_A}}        \def\RRa{\Ro{{\cal R}^-_A}}
\def\Rb{\Ro{{\cal R}_B}}        \def\RRb{\Ro{{\cal R}^-_B}}
\def\Rh{\Ro{{\cal R}_H}}
\def\id{\kr
\begin{picture}(0,2)\put(0,0){\line(0,1)2}\end{picture}}
\def\obj##1{\settowidth{\textwd}{$##1$}%
\raise .2\unitlength\hbox{\kern -.5\textwd $##1$ \kern -.5\textwd \krrr}}
\def\Obj##1{\settowidth{\textwd}{$##1$}%
\raise 1.1\unitlength\hbox{\kern -1\textwd $##1$}}
\def\hhbox##1{\hbox{\krrrr
\def\coev{\kr
\begin{picture}(1,1)\put(.5,0){\oval(1,1)[t]}\end{picture}}
\def\ev{\kr
\begin{picture}(1,1)\put(.5,1){\oval(1,1)[b]}\end{picture}}
\def\ld{\kr
\begin{picture}(1,1)
\put(1,0){\oval(2,2)[lt]}\put(1,0){\line(0,1)1}
\end{picture}}
\def\Ld{\kr
\begin{picture}(2,1)
\put(2,0){\oval(4,2)[lt]}\put(2,0){\line(0,1)1}
\end{picture}}
\def\rd{\kr
\begin{picture}(1,1)
\put(0,0){\oval(2,2)[rt]}\put(0,0){\line(0,1)1}
\end{picture}}
\def\Rd{\kr
\begin{picture}(2,1)
\put(0,0){\oval(4,2)[rt]}\put(0,0){\line(0,1)1}
\end{picture}}
\def\cd{\kr
\begin{picture}(1,1)
\put(.5,0){\oval(1,1)[ct]}\put(.5,.5){\line(0,1){.5}}
\end{picture}}
\def\lu{\kr
\begin{picture}(1,1)
\put(1,1){\oval(2,2)[lb]}\put(1,0){\line(0,1)1}
\end{picture}}
\def\Lu{\kr
\begin{picture}(2,1)
\put(2,1){\oval(4,2)[lb]}\put(2,0){\line(0,1)1}
\end{picture}}
\def\cu{\kr
\begin{picture}(1,1)
\put(.5,1){\oval(1,1)[cb]}\put(.5,0){\line(0,1){.5}}
\end{picture}}
\def\ru{\kr
\begin{picture}(1,1)
\put(0,1){\oval(2,2)[rb]}\put(0,0){\line(0,1)1}
\end{picture}}
\def\Ru{\kr
\begin{picture}(2,1)
\put(0,1){\oval(4,2)[rb]}\put(0,0){\line(0,1)1}
\end{picture}}
\def\hru{\kr
\begin{picture}(.5,1)
\put(0,1){\oval(1,1)[rb]}\put(0,0){\line(0,1)1}
\end{picture}}
\def\hlu{\kr
\begin{picture}(.5,1)
\put(.5,1){\oval(1,1)[lb]}\put(.5,0){\line(0,1)1}
\end{picture}}
\def\hrd{\kr
\begin{picture}(.5,1)
\put(0,0){\oval(1,1)[rt]}\put(0,0){\line(0,1)1}
\end{picture}}
\def\hld{\kr
\begin{picture}(.5,1)
\put(.5,0){\oval(1,1)[lt]}\put(.5,0){\line(0,1)1}
\end{picture}}
\def\id{\kr
\begin{picture}(0,1)\put(0,0){\line(0,1)1}\end{picture}}
\def\d{\kr
\begin{picture}(.5,1)\put(0,1){\line(1,-2){0.5}}\end{picture}}
\def\dd{\kr
\begin{picture}(.5,1)\put(0,0){\line(1,2){0.5}}\end{picture}}
##1}}#1}\normalbaselines}}
\def\object#1{\settowidth{\textwd}{$#1$}%
                        \hbox{%
                        \kern -.5\textwd $#1$ \kern -.5\textwd}}
\def\map#1#2#3{\vcenter{\hbox{$#2\;$}}
                     \vcenter{\settowidth{\textwd}{$#1$}
	                      \hbox{\kern -.5\textwd $#1$ \kern -.5\textwd}
			      \hbox{\begin{picture}(0,2)
                                          \put(0,2){\vector(0,-1)2}
                                    \end{picture}}
                              \settowidth{\textwd}{$#3$}
	                      \hbox{\kern -.5\textwd $#3$ \kern -.5\textwd}}}

\begin{document}
\begin{center}
\Large { \bf On braided FRT-construction}
\end{center}
\bigskip
{\bf Yuri BESPALOV }
\footnote{The research described in this paper was made possible
          by Grant No U4J200 from the International Science
Foundation.}
\bigskip \\ {\it Bogolyubov Institute for Theoretical Physics }
\par\noindent
{\it Metrologichna str., 14-b \ \ Kiev 143, 252143 Ukraine}
\par\noindent
{\it E-mail: mmtpitp@gluk.apc.org}
\par\medskip\noindent
{\it Received: September, 1995}
\bigskip
\begin{abstract}
\small
Fully braided analog of Faddeev-Reshetikhin-Takhtajan construction of
quasitriangular bialgebra $A(X,R)$ is proposed.
For given pairing $C$ factor-algebra $A(X,R;C)$ is a dual quantum
braided group.
Corresponding inhomogeneous quantum group is obtained as a result of
generalized bosonization.
Construction of first order bicovariant differential calculus
is proposed.
\end{abstract}

\section{Introduction and preliminaries}

Hopf algebras in braided categories (braided groups) have been extensively
studied over the last few years and play an important role in
$q$-deformed physics and mathematics \cite{Majid8},\cite{Majid10}.
Examples, applications and the basic theory of braided groups have been
introduced and developed by Majid.
Some similar concepts arise independently in works of Lyubashenko
inspired by results on conformal field theory.
Crossed modules over braided Hopf algebras were introduced and studied
in \cite{Bespalov2} and provide a useful technique for investigation
of braided Hopf algebras.
In particular, crossed product of braided Hopf algebras and generalized
bosonization for quantum braided groups are defined in \cite{Bespalov2}.
The theory of Hopf bimodules in braided categories is developed in
\cite{BD1}, on grounds of \cite{Bespalov2}.
Application of this theory is an analog of Woronowicz construction of
(bicovariant) differential calculi \cite{Wor} developed in \cite{BD} for
the case of braided Hopf algebras and quantum braided groups.
Quantum braided group defined by Majid \cite{Majid7,Majid6}
is a natural generalization of Drinfel'd's concept of (ordinary) quantum
group (quasitriangular Hopf algebra) \cite{Drinfel'd1}.
Basic examples of coquasitriangular bialgebras $A(R)$ are obtained as a
result of Faddeev-Reshetikhin-Takhtajan construction \cite{FRT}
applied to an arbitrary $R$-matrix.
Analog of FRT-construction for anyonic quantum groups is described in
\cite{MR}.
Majid proposed another construction of braided bialgebra $B(R)$ which
can be obtained as a transmutation \cite{Majid7} of $A(R)$.
Algebra $A(R,Z)$ defined in \cite{H} generalizes both  $A(R)$ and $B(R)$.
In this paper we describe a fully braided analog of
FRT-construction of quasitriangular bialgebra $A(X,R)$, where $X$ is an
object of an Abelian braided monoidal category $\C$ and
$R:X\otimes X\rightarrow X\otimes X$ solution of the braid equation.
This construction covers all mentioned above and can be considered as
a coordinate-free version of \cite{H}.
For a given pairing $C$ we define a factor-algebra $A(X,R;C)$ which is
a dual quantum braided group.
This is an analog of construction of quantum simple Lie groups of type
$B,C,D$ in \cite{FRT}.
Majid's definition of braided vectors $\V(R)$ is simply reformulated to our
more abstract setting.
In particular, our algebra $\V(X,R)$ is also a quantum braided
group in the category of comodules over $A(X,R)$.  Quantized analogs of
inhomogeneous linear groups are studied in many papers (see \cite{AC,D}
and references therein).  Generalized bosonization construction
\cite{Bespalov2} allows us to define quantum braided group
$A(X,R)\ltimes\V(X,R)$.  We propose construction of a first order
bicovariant differential calculus on dual quantum braided group $A$
related with any comodule $X$ over $A$.  In our special case
$A=A(X,R)$ this is a
generalization of construction \cite{J}.

In the rest of this part we give necessary preliminary results.
The main results of the paper are presented in the second part.

\subsec{}
We will suppose that
${\cal C}$
is {\em an Abelian}
and {\em braided (monoidal) category}
with tensor product  $\otimes$, unit object $\underline 1$
and braiding $\Psi$
(without loss of generality by Mac Lane's coherence theorem
we will assume that underlying monoidal category is strict, i.e. the functors
$\_\otimes (\_\otimes\_)$ and $(\_\otimes\_)\otimes\_$ coincide and
$\underline 1\otimes X=X=X\otimes\underline 1$).
Compatibility conditions between tensor product and Abelian structure
are the following \cite{BD}:
functors $(-)\otimes X$ and $X\otimes (-)$ are right exact for any object $X$
(this assumption is true if the category is closed);
for any epimorphisms $X_i\buildrel{f_i}\over\rightarrow Y_i,\enspace i=1,2$
the diagram
\begin{equation}
          \label{puth-out}
          \matrix{ X_1\otimes X_2 &
          \buildrel{X_1\otimes f_2}\over\longrightarrow  &
	  X_1\otimes Y_2 \cr
	  {}^{f_1\otimes X_2}\downarrow &&
	  \downarrow^{f_1\otimes Y_2} \cr
	  Y_1\otimes X_2 &
          \buildrel{Y_1\otimes f_2}\over\longrightarrow  &
	  Y_1\otimes Y_2}
\end{equation}
is push-out (the right-down part is a colimit of the left-up part).
In this case there exists well-behaviored constructions of
factor-algebra (coalgebra, bialgebra, Hopf algebra) by ideal (coideal,
biideal, Hopf ideal).
One can define an algebra by generator and relations.
We means under 'the ideal generated by relations
$f_1=f_2:X\rightarrow A\,$' the subobject
${\rm Im}\left(\mu\circ(\mu\otimes A)\circ
(A\otimes(f_1-f_2)\otimes A)\right)$ of algebra $A$.

\subsec{}
We will work with {\em graded and filtered algebras in $\cal C$}.
A ($\NN$-)graded algebra $A$ means a collection of objects
$A_k,\;k\in{\cal C}$, multiplications
$m_{i,j}:A_i\otimes A_j\rightarrow A_{i+j}$ satisfying associativity
conditions and unit $\eta:\underline 1\rightarrow A_0$.
A ($\NN$-)filtered algebra $A$ means a collection of objects
$A_{(k)},\;k\in{\cal C}$, such that $A_{(i)}$ is subobject of
$A_{(j)}$ if $i<j$, multiplications $m_{(i),(j)}:A_{(i)}\otimes
A_{(j)}\rightarrow A_{(i+j)}$ satisfying conditions of associativity
and compatibility with restrictions on subobjects, and unit
$\eta:\underline 1\rightarrow A_{(0)}$.
For any graded algebra $\{A_i\}$ the collection
$\{A_{(k)}:=\oplus_{i=0}^kA_i\}$ with natural multiplications is a
filtered algebra.
As shown in \cite{BD} graded or filtered algebra can be considered as
a usual algebra in a certain category of 'graded spaces' i.e. functors
from a certain category to category $\cal C$. This category of 'graded
spaces' is again an Abelian braided monoidal category.
Similarly graded coalgebras, bialgebras, Hopf algebras can be defined.
We will say briefly that graded (filtered) algebra lives in a category
$\C$ if its components $A_n$ ($A_{(n)}$) live in $\C$.
See \cite{BD} about more details.

\subsec{}
We actively use diagrammatic calculus in braided categories
\cite{Majid6,Majid8} (see \cite{Bespalov2} about our slight
modifications).
Morphisms $\Psi$ and $\Psi^{-1}$ are represented by
under and over crossing and algebraic information 'flows' along braids
and tangles according to functoriality and the coherence theorem for
braided categories \cite{JS}:
\begin{equation}
\Psi=\enspace
\vvbox{\hbox{\hx}}
\qquad\quad
\Psi^{-1}=\enspace
\vvbox{\hbox{\hxx}}
\qquad\qquad
\vvbox{\hbox{\O{f}\step\id}
       \hbox{\hx}}
\enspace =\enspace
\vvbox{\hbox{\hx}
       \hbox{\id\step\O{f}}}
\qquad\quad
\vvbox{\hbox{\id\step\O{f}}
       \hbox{\hx}}
\enspace =\enspace
\vvbox{\hbox{\hx}
       \hbox{\O{f}\step\id}}
\label{Psi}
\end{equation}

\begin{figure}
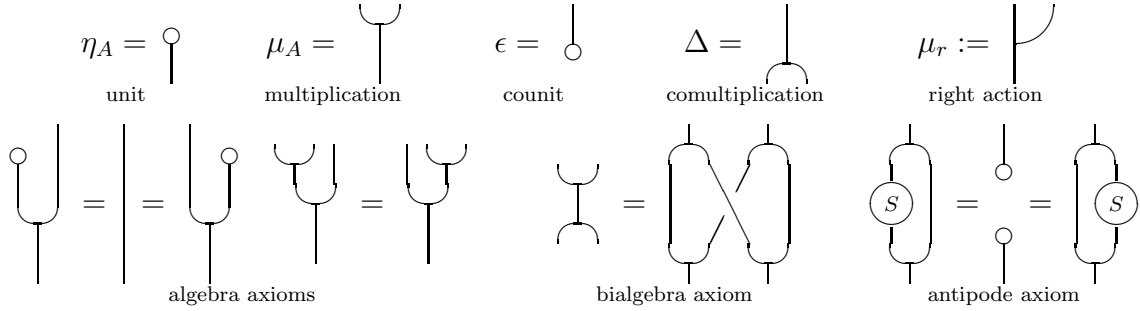

$$
\matrix{
\eta_A={}\enspace
\vvbox{\hbox{\unit}}
&\quad&
\mu_A={}\enspace{}
\vvbox{\hbox{\hcu}}
&\quad&
\epsilon ={}\enspace{}\;
\vvbox{\hbox{\counit}}
&\quad&
\Delta ={}\enspace{}
\vvbox{\hbox{\hcd}}
&\quad&
\mu_r:={}\enspace
\vvbox{\hbox{\ru}}
\cr
\hbox{\scriptsize unit}
&&
\hbox{\scriptsize multiplication}
&&
\hbox{\scriptsize counit}
&&
\hbox{\scriptsize comultiplication}
&&
\hbox{\scriptsize right action}
}
$$
$$
\matrix{
\vvbox{\hbox{\unit\step\id}\hbox{\hcu}}
\enspace ={}\;
\vvbox{\hbox{\id}\hbox{\id}}
\;={}\enspace{}
\vvbox{\hbox{\id\step\unit}\hbox{\hcu}}
\quad\enspace
\vvbox{\hhbox{\cu\hstep\id}
       \hbox{\hstep\hcu}}
\enspace ={}\enspace{}
\vvbox{\hhbox{\id\hstep\cu}
       \hbox{\hcu}}
&\quad&
\vvbox{\hhbox{\cu}
       \hhbox{\cd}}
\enspace ={}\enspace{}
\vvbox{\hhbox{\cd\step\cd}
       \hbox{\id\step\hx\step\id}
       \hhbox{\cu\step\cu}}
&\quad\;&
  \vvbox{\hhbox{\cd}
       \hbox{\S\step\id}
       \hhbox{\cu}}
\enspace ={}\enspace{}
\vvbox{\hbox{\counit}
       \hbox{\unit}}
\enspace ={}\enspace{}
\vvbox{\hhbox{\cd}
       \hbox{\id\step\S}
       \hhbox{\cu}}
\cr
\hbox{\scriptsize algebra axioms}
&&
\hbox{\scriptsize bialgebra axiom}
&&
\hbox{\scriptsize antipode axiom}
}
$$
\caption{ The basic algebraic structures in a braided category }
\label{Fig-Main}
\end{figure}

Fig.\ref{Fig-Main} explains our notations:
{\em An algebra} in a monoidal category $\cal C$ is an object $A$
equipped with unit $\eta=\eta_A:\, \underline 1\rightarrow A$
and multiplication $\mu=\mu_A:\, A\otimes A\rightarrow A$
obeying the axioms on Fig.\ref{Fig-Main}.
{\em A coalgebra} is object $C$ equipped with
counit $\epsilon=\epsilon_A:\,C\rightarrow\underline 1$ and
comultiplication $\Delta=\Delta_A:\,A\rightarrow A\otimes A$
obeying the axioms of algebra turned upside-down
Finally, \cite{M2},\cite{Majid7}
{\em a bialgebra $A$ in a braided category $\cal C$} is
an object in $\cal C$ equipped with algebra and coalgebra structures
obeying the compatibility axiom on Fig.\ref{Fig-Main} which means
that $\Delta_A$ is  an algebra homomorphism.
{\em A Hopf algebra $A$ in a braided category} $\cal C$
({\em braided group\/} or {\em braided Hopf algebra}) is
a bialgebra in $\cal C$ with antipode $S:\,A\rightarrow A$ which is
convolution-inverse to identical map
(the last identity on Fig.\ref{Fig-Main}).
Axioms for (co-)module $X$ over a (co-)algebra $A$
are obtained by "polarization" of the (co-)algebra axioms.

If $\C$ is a braided category we will denote by $\overline\C$
the same category with the same tensor product and with
inverse braiding $\Psi^{-1}$.
For any algebra (resp. coalgebra) $A$ in $\cal C$ we will always consider
{\em the opposite algebra\/}
$(A^{\rm op},\mu_{A^{\rm op}}:=\mu_A\circ\Psi^{-1})$
(resp. {\em the opposite coalgebra\/}
$(A_{\rm op},\Delta_{A_{\rm op}}:=\Psi^{-1}\circ\Delta_A$)
as an object of the category $\overline{\cal C}$.
In particular, $(A^{\rm op})^{\rm op}=A$.
If $A$ is a bialgebra in $\cal C$ then $A^{\rm op}$ and $A_{\rm op}$ are
bialgebras in $\overline{\cal C}$ (cf. \cite{Majid8}).
Antipode $S^-$ for $A^{\rm op}$ (or, the same, for $A_{\rm op}$) is called
{\em skew antipode} and equals $S^{-1}$ if both $S$ and $S^-$
exist. Majid \cite{Majid8} derived from Hopf algebra axioms that
antipode $S_A$ is a bialgebra morphism $(A^{\rm op})_{\rm op}\rightarrow A$
(or $A\rightarrow (A_{\rm op})^{\rm op}$) in $\cal C$.

\begin{figure}
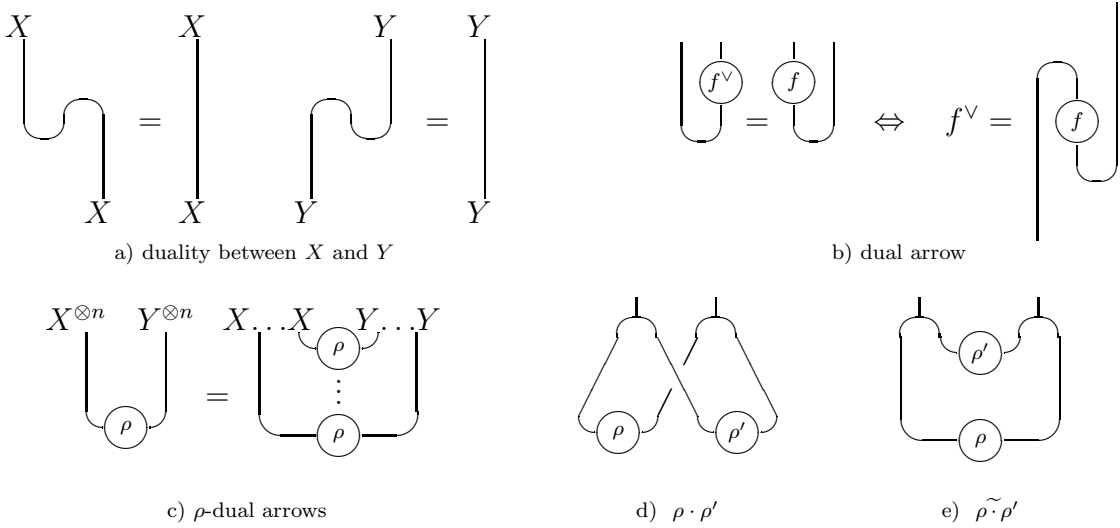

$$
\matrix{
\matrix{\object{X}\Step\cr
	\vvbox{\hbox{\id\step\hcoev}
	       \hbox{\hev\step\id}}\cr
	\Step\object{X}}
\enspace =\enspace
\matrix{\object{X}\cr
        \vvbox{\hbox{\id}\hbox{\id}}\cr
        \object{X}}
\qquad\quad
\matrix{\Step\object{Y}\cr
	\vvbox{\hbox{\hcoev\step\id}
	       \hbox{\id\step\hev}}\cr
        \object{Y}\Step}
\enspace =\enspace
\matrix{\object{Y}\cr
        \vvbox{\hbox{\id}\hbox{\id}}\cr
        \object{Y}}
&\qquad\qquad&
\vvbox{\hbox{\id\step\O{f^\vee}}
       \hbox{\hev}}
\enspace =\enspace
\vvbox{\hbox{\O{f}\step\id}
       \hbox{\hev}}
\quad\Leftrightarrow\quad
f^\vee =\enspace
\vvbox{\hbox{\hcoev\step\id}
       \hbox{\id\step\O{f}\step\id}
       \hbox{\id\step\hev}}
\cr
\hbox{\scriptsize a) duality between $X$ and $Y$}
&&
\hbox{\scriptsize b) dual arrow}
}
$$
$$
\matrix{
\matrix{\object{X^{\otimes n}}\Step\object{Y^{\otimes n}}\cr
\vvbox{\hbox{\id\Step\id}
       \hbox{\coro{\rho}}}}
\enspace=\enspace
\matrix{\object{X\!\dots\! X}\step\Step\object{Y\!\dots\! Y}\cr
\vvbox{\hbox{\id\step\coro{\rho}\hhstep\hhstep\obj{\vdots}\step\step\id}
       \hbox{\coRo{\rho}}}}
&\qquad\quad&
\vvbox{\hhbox{\step\cd\step\cd\step}
       \hbox{\dd\step\hx\step\d}
       \hbox{\coro{\rho}\step\coro{\rho^\prime}}}
&\qquad&
\vvbox{\hhbox{\cd\Step\cd}
       \hbox{\id\step\coro{\rho^\prime}\step\id}
       \hbox{\coRo{\rho}}}
\cr
\hbox{\scriptsize c) $\rho$-dual arrows}
&&
\hbox{\scriptsize d) \ $\rho\cdot\rho^\prime$}
&&
\hbox{\scriptsize e) \ $\rho\,\widetilde\cdot\,\rho^\prime$}
}
$$
\caption{Duals and pairings.}
\label{Fig-Pairing}
\end{figure}

\subsec{}
For objects $X,Y$ of a monoidal category $\cal C$ we will call any
morphism
\begin{equation}
\cup=\cup^{X,Y} :\,X\otimes Y\rightarrow\underline 1\qquad
(\;{\rm resp.}\enspace
\cap=\cap_{Y,X} :\,\underline 1\rightarrow Y\otimes  X\;)
\label{Equation-Pairing}
\end{equation}
{\em a pairing between $X,Y$} (resp. {\em copairing between $Y,X$}).
{\em Duality between $X$ and $Y$} is both pairing and copairing
(\ref{Equation-Pairing}) obeying the identities on
Fig.\ref{Fig-Pairing}a.
In this case $X$ is called {\em left dual} to $Y$ (resp. $Y$ is called
{\em right dual} to $X$) and we will write $X={}^\vee Y, Y=X^\vee$.
Dual arrow $f^\vee$ is defined by one of the two equivalent
conditions on Fig.\ref{Fig-Pairing}b.  In this way a braided monoidal
functor $(\_)^\vee:\,{\cal C}\rightarrow{\cal C}^{\rm op}_{\rm op}$ can
be defined if $X^\vee$ exists for each $X\in{\rm Obj}({\cal C})$.
Without loss of generality by coherence theorem we shall assume that
$(\_)^\vee$ is a strict monoidal functor:
$(X\otimes Y)^\vee=Y^\vee\otimes X^\vee$,
$(f\otimes g)^\vee=g^\vee\otimes f^\vee$.
Pairing $\rho$ between $X$ and $Y$ extends to pairing  between
$X^{\otimes n}$ and $Y^{\otimes n}$ defined by
the diagram on Fig.\ref{Fig-Pairing}c.
We say that arrows
$f:\,X^{\otimes m}\rightarrow X^{\otimes n}$ and $g:\,Y^{\otimes
n}\rightarrow Y^{\otimes m}$ are $\rho$-{\em dual} if
$\rho\circ(f\otimes Y^{\otimes n})=\rho\circ(X^{\otimes n}\otimes g)\,.$

Let $A$ and $H$ be bialgebras in braided category $\cal C$.
Morphism $\rho:\,A\otimes H\rightarrow\underline 1$ is called
{\em a bialgebra pairing}
if algebra (resp. coalgebra) structure on $A$ and
coalgebra (resp. algebra) structure on $H$ are $\rho$-dual.
Convolution product '$\cdot$' and 'the second' product
'$\widetilde\cdot$' for $\rho,\rho^\prime\in{\rm Hom}_{\cal C}(X\otimes
Y,\underline 1)$ are defined on Fig.\ref{Fig-Pairing}d,e.
We denote by
$\rho^-$, $\rho^{\sim{}}$ corresponding inverse to $\rho$.  Let
$\overline\rho:=\rho^-\circ\Psi^{-1}$.  If $A$ or $H$ has (skew)
antipode then $\rho^{\sim{}}$ (resp. $\rho^-$) exists and
\begin{equation}
\rho\circ(S_A\otimes H)=\rho^{\sim{}}=\rho\circ(A\otimes S_H)
\qquad\enspace
   \rho\circ(S_A^-\otimes H)=\rho^-=\rho\circ(A\otimes S_H^-)
\end{equation}
If $\rho^-$ or $\rho^{\sim{}}$ exists then $\rho$-duality between
multiplications and comultiplications implies $\rho$-duality between
units and counits.
If $(A,H,\rho)$ is bialgebra pairing in $\cal C$ then
$(A_{\rm op},H_{\rm op},\rho^-)$, $(A^{\rm op},H^{\rm op},\rho^{\sim{}})$,
$(H^{\rm op},A^{\rm op},\overline\rho)$ are bialgebra pairings in
$\overline{\cal C}$.

\subsec{}
Quantum braided groups in a braided category were introduced in \cite{Majid7}
and basic theory was developed there.
The following are input-output reversed variant of definitions from
\cite{Majid7} in a slightly modified form
\cite{Bespalov2} suitable for our use.

\begin{figure}
$$
\matrix{
\vvbox{\hbox{\cd\step\cd}
       \hbox{\id\Step\hx\Step\id}
       \hbox{\cu\hstep\obj{\overline\mu^{\rm op}}\hstep\coro{\rho}}}
\enspace=\enspace
\vvbox{\hbox{\cd\step\cd}
       \hbox{\id\Step\hx\Step\id}
       \hbox{\coro{\rho}\step\cu}  }\cr
\hbox{\scriptsize a) The axiom for a dual quantum braided group
                 $(A,\overline A,{\rho})$. }
}
$$
$$
\matrix{
\matrix{
        \object{A}\step\object{X}\step
  \cr
\vvbox{
       \hbox{\id\step\rd}
       \hbox{\hx\step\id}
  \hhbox{\krl\id\step\cu\hstep\obj{\mu}}
  }\cr
  \object{X}\step\hstep\object{A}\hstep}
\enspace=\enspace
\matrix{
        \object{A}\step\object{X}\step
  \cr
\vvbox{
       \hbox{\id\step\rd}
       \hbox{\hxx\step\id}
  \hhbox{\krl\id\step\cu\hstep\obj{\overline\mu^{\rm op}}}
  }\cr
  \object{X}\step\hstep\object{A}\hstep
  }
&\qquad\qquad&
\Psi=\enspace
\vvbox{
       \hbox{\id\step\rd}
       \hbox{\hx\step\d}
       \hhbox{\krl\id\step\hrd\step\hstep\d}
  \hbox{\id\step\id\hstep\coro{\rho}}
  }
\qquad\quad
\Psi^{-1}=\enspace
\vvbox{
       \hbox{\id\step\rd}
       \hbox{\hxx\step\d}
       \hhbox{\krl\id\step\hrd\step\hstep\d}
  \hbox{\id\step\id\hstep\coro{\overline{\rho}}}
  }
\cr
\hbox{\scriptsize b) the condition on comodules from ${\cal C}^\cO{A}$}
&&
\hbox{\scriptsize c) braiding in ${\cal C}^\cO{A}$}
}
$$
\caption{}
\label{Fig-QBG}
\end{figure}

{\em A coquasitriangular bialgebra} in a braided category $\cal C$ is a
  pair of bialgebras $A$ in $\cal C$ and $\overline A$ in
$\overline{\cal C}$ with the same underlying coalgebra ($\mu$ and
$\overline\mu$ are multiplications in $A$ and $\overline A$
 respectively), and convolution invertible bialgebra pairing ({\em
coquasitriangular structure})
${\rho}:\, \overline A^{\rm op}\otimes A
\rightarrow\underline 1$,
satisfying the condition on
Fig.\ref{Fig-QBG}a.  (It follows directly from the definition that
units for $A$ and for $\overline A$ are the same.) {\em A dual quantum
braided group} or {\em a coquasitriangular Hopf algebra} in $\cal C$ is
a coquasitriangular bialgebra such that $A$ and $\overline A$ have
antipodes $S$ and $\overline S$ respectively.  (In this case
${\rho}^-= \rho\circ(\overline S\otimes A)$
         and
$\rho^{\sim{}}=\rho\circ(A\otimes S)$.)
         \par
In particular, for any bialgebra (braided group) $A$ the pair $(A,A^{\rm
 op})$  is a coquasitriangular bialgebra (dual quantum braided group)
with the trivial coquasitriangular structure
$\rho=\epsilon\otimes\epsilon$.

Category ${\cal C}^\cO{A,\overline A}$ is
a full subcategory of the category ${\cal C}^A$ of right comodules with
         objects $X$ satisfying the first identity on
         Fig.\ref{Fig-QBG}b.  ${\cal C}^\cO{A,\overline A}$ is a
         monoidal subcategory of ${\cal C}^A$ and braided with $\Psi$
and $\Psi^{-1}$ shown on Fig.\ref{Fig-QBG}c.
We use a brief notation ${\cal C}^\cO{A}$ for
${\cal C}^\cO{A,A^{\rm op}}$.

\section{On braided FRT-construction}

\subsec{}
Canonical epimorphism $B_n\rightarrow S_n$ of the braid group
into the permutation group
admits a section $S_n\buildrel{\widehat{}}\over\rightarrow B_n$
identical on generators and unqueenly defined by condition that
$\widehat{\sigma_1\sigma_2}=\widehat{\sigma_1}\widehat{\sigma_2}$ if
$\ell(\sigma_1\sigma_2)=\ell(\sigma_1)+\ell(\sigma_2)$ where
$\ell(\sigma)$ is the length (of the minimal decomposition) of $\sigma$.
For any object $X$ obvious action of braid group $B_n$ on
$X^{\otimes n}$ is defined.
We will use the same notation $\widehat\sigma$ for the image of
the braid $\widehat\sigma\in B_n$ in
${\rm End}_{\cal C}(X^{\otimes n})$.
For $k=1,\dots,n$ let us denote by $S^k_n\subset S_n$ subset of
${n!}\over{k!(n-k)!}$ shuffle
permutations which preserves order of any two elements $i$ and $j$ if
$i,j\le k$ or $i,j>k$.
Majid in \cite{Majid11} defines braided binomial coefficient
as a sum of ${n!}\over{k!(n-k)!}$ braids in
${\rm End}_{\cal C}(X^{\otimes n})$ and in particular, braided
factorial as a sum of $n!$ braids:
\begin{equation}
\left[{n\atop k}; X\right] :=
\sum_{\sigma^{-1}\in S_n^k}\widehat\sigma, \qquad
[n;X]!:=\sum_{\sigma\in S_n}\widehat\sigma.
\end{equation}

\subsec{}
\label{tensor-Hopf}
For any object $X$ of a braided category $\cal C$ the tensor algebra
$\T{X}=\{ X^{\otimes n}\}_{n\in{\bf N}}$ is
a graded Hopf algebra
with the tensor product as multiplication, comultiplication
\begin{equation}
\Delta_{m,n}:=\left[{{m+n}\atop {m}};X\right]:\;
  X^{\otimes (m+n)}\rightarrow X^{\otimes m}\otimes X^{\otimes n}
\end{equation}
and antipode
\begin{equation}
S\vert_{X^{\otimes n}}:=
  (-1)^n\circ\widehat{\rho_n}:\;
X^{\otimes n}\rightarrow X^{\otimes n}.
\end{equation}
Where $S_n\ni\rho_n:(1,2,\dots ,n)\mapsto(n,n-1,\dots ,1)$ and
$\widehat{\rho_n}$ is a Garside element of $B_n$.
The bialgebra axiom turns into Newton-Majid
binomial formula \cite{Majid11}:
\begin{equation}
(\underline 1\otimes x+x\otimes\underline 1)^n=
(\Delta\circ x)^n=
\sum^n_{k=0}\left[{n\atop k};X\right]\circ
(x^{\otimes k}\otimes x^{\otimes (n-k)}) \qquad
\end{equation}
for any $x:\; Z\rightarrow X$.

One can define two graded ideals in $\T{X}$:
$I=\{ I_n\subset X^{\otimes n}\}$ is an ideal
generated by its 'quadratic part'
$I_2:={\rm ker}\,[2;X]!={\rm ker}(\Psi_{X,X}+{\rm id})$.
And let $I^\bullet=\{ I^\bullet_n:={\rm ker}\,[n; X]\}$.
\par
$I$ is non zero iff $-1$ is an eigenvalue of $\Psi_{X,X}$.
If we suppose that the multiplicity of this eigenvalue is $1$, i. e.
we can choose a minimal polynomial of $\Psi_{X,X}$ in
the form $p(t)=p_{-1}(t)(t+1)$ with $p_{-1}(-1)=1$,
then $P_{-1}:=p_{-1}(\Psi_{X,X})$ is an idempotent, and
in this case $I_2={\rm ker}\,[n;X]!={\rm im}\,P_{-1}$.
\par
It's easy to see that $I\subset I^\bullet$
The following example from \cite{BD} show that, in general, the ideal
$I^\bullet$ has 'generators' of power more than $2$.

\begin{example}
Let $\T{X}$ be an algebra generated by the one dimensional vector
space $X=kx$ over a field $k$ with the braiding
$\Psi(x\otimes x)=q(x\otimes x)$,
$q\in k$.
In this case braided integers are 'ordinary' $q$-integers:
$[n]_q:=1+q+\dots +q^{n-1}\,.$
And for $q$ a primitive root of $1$ of order $n>2$:
$I=\emptyset$ but $I^\bullet=(x^{\otimes n})$.
\end{example}

\begin{proposition}
Both $I$ and $I^\bullet$ are Hopf ideals in $\T{X}$.
We denote by $\V(X)$ and $\V^\bullet(X)$ corresponding factor-algebras
\end{proposition}

For the special case of the category $\C$ built from an arbitrary $R$-matrix
$\V(X)$ is an algebra of functions on 'quantum vector space'.
Majid discovered a Hopf algebra structure on this object (cf. \cite{Majid11}
and references therein).

\subsec{}
It is well known that any solution
$R=R_{X,X}:X\otimes X\rightarrow X\otimes X\,,\;X\in{\rm Obj}(X)$
of the braid equation
\begin{equation}
(X\otimes R)\circ (R\otimes X)\circ (X\otimes R)=
(R\otimes X)\circ (X\otimes R)\circ (R\otimes X)
\end{equation}
with a certain invertibility conditions defines a braided structure
on the monoidal subcategory of $\cal C$ generated by object $X$ and its dual
${}^\vee\!X$ as described in what follows.
Morphisms $R_{X^{\otimes m},X^{\otimes n}}$ are uniquely defined by the
hexagon identities:
\begin{equation}
\label{Eq-Hex-id}
{R}_{Y\!\otimes\!Y^\prime,Z}=
({R}_{Y,Z}\otimes Y^\prime)\circ (Y \otimes{R}_{Y^\prime,Z})\,,
\quad
{R}_{Y,Z\!\otimes\!Z^\prime}=
(Z\otimes{R}_{Y,Z^\prime})\circ(R_{Y,Z}\otimes Z^\prime)\,.
\end{equation}
where $Y,Y^\prime,Z,Z^\prime$ are powers of $X$.
And let
$R_{{}^\vee\!X^{\otimes m},{}^\vee\!X^{\otimes n}}:=
 {}^\vee\left(R_{X^{\otimes m},X^{\otimes n}}\right)$.
We also suppose that there exists
$R_{X,{}^\vee\!X}:X\otimes{}^\vee\!X\rightarrow {}^\vee\!X\otimes X$
inverse to
$(\cup\otimes X\otimes{}^\vee\!X)\circ
 ({}^\vee\!X\otimes R_{X,X}\otimes {}^\vee\!X)\circ
 ({}^\vee\!X\otimes X\otimes \cap)$
Let $R_{X^{\otimes m},{}^\vee\!X^{\otimes n}}$ be uniquely defined by
the hexagon identities (\ref{Eq-Hex-id}) and
$R_{{}^\vee\!X^{\otimes m},X^{\otimes n}}$ be defined in dual way.
$\C(X,{}^\vee\!X;R)$ is a subcategory of $\C$ whose objects are tensor
products of $X$ and ${}^\vee\!X$ and morphism $f:Y\rightarrow Z$
are those from $\C$ which 'flow' along the braids labeled by
$R_{X,\_}$ and $R_{\_,X}$, i.e.
\begin{equation}
R_{X,Z}\circ(X\otimes f)=(f\otimes X)\circ R_{X,Y}\,,\quad
R_{Z,X}\circ(f\otimes X)=(X\otimes f)\circ R_{Y,X}\,.
\end{equation}

The analog $A(X,R)$ of FRT-bialgebra \cite{FRT} can be obtained as a result
of a some type of reconstruction for the monoidal functor
$\C(X,{}^\vee\!X;R)\rightarrow\C$.

\subsec{}
As the first step the following lemmas allow us to define a bialgebra
$A(X)$.

\begin{lemma}
Let $X$ be an object of $\cal C$ with (left) dual ${}^\vee\!X$.
Then ${}^\vee\!X\!\otimes\! X$ can be equipped with a coalgebra structure
\begin{equation}
\Delta_{{}^\vee\!X\!\otimes\! X}:=
{}^\vee\!X\otimes\cap_{X,{}^\vee\!X}\otimes X\,,\quad
\epsilon_{{}^\vee\!X\!\otimes\! X}:=\cup^{{}^\vee\!X,X}\,.
\end{equation}
$X$ (resp. ${}^\vee\!X$) becomes a right (resp. left) comodule over
${}^\vee\!X\!\otimes\! X$ with coaction
\begin{equation}
\label{Delta-r}
\Delta^X_r:=\cap_{X,{}^\vee\!X}\otimes X\,,\quad
\Delta^{{}^\vee\!X}_\ell:={}^\vee\!X\otimes\cap_{X,{}^\vee\!X}\,.
\end{equation}
\end{lemma}

\begin{lemma}
Let $(A,\Delta_A)$, $(B,\Delta_B)$ be coalgebras in $\cal C$ and
$(X,\Delta^X_r)$, $(Y,\Delta^Y_r)$ right comodules
over $A$ and $B$ respectively. Then $A\!\otimes\! B$ is a coalgebra
with comultiplication
\begin{equation}
\Delta_{A\otimes B}:=
(A\otimes\Psi_{A,B}\otimes B)\circ(\Delta_A\otimes\Delta_B)
\end{equation}
and
$X\!\otimes\! Y$ is a right $A\!\otimes\! B$-comodule
with coaction
\begin{equation}
\Delta^{X\otimes Y}_r:=
(X\otimes\Psi_{A,Y}\otimes B)\circ(\Delta^X_r\otimes\Delta^Y_r)
\end{equation}
\end{lemma}

So with any two objects $X$ and $Y$ which have left duals we can connect
the following coalgebras: tensor product of two coalgebras
$({}^\vee\!X\otimes X)\otimes({}^\vee Y\otimes Y)$
and coalgebra $({}^\vee Y\otimes{}^\vee\!X)\otimes (X\otimes Y)$
related with the object $X\otimes Y$.

\begin{lemma}
Morphism
\begin{equation}
\mu_{{}^\vee\!X\otimes X,{}^\vee\!Y\otimes Y}:=
\Psi_{{}^\vee\!X\otimes X,{}^\vee\!Y}\otimes Y:\,
 ({}^\vee\!X\otimes X) \otimes({}^\vee\!Y\otimes Y)
 \rightarrow
 ({}^\vee\!Y\otimes{}^\vee\!X)\otimes (X\otimes Y)
\end{equation}
is coalgebra isomorphism and interlaces
coactions of these coalgebras on $X\otimes Y$.
\par
For objects
$X,Y,Z$ with left dual the following {\em associativity condition} is
true:
\begin{equation}
\mu_{X\otimes Y,Z}\circ (\mu_{X,Y}\otimes Z)=
\mu_{X,Y\otimes Z}\circ (X\otimes \mu_{Y,Z})
\end{equation}
\end{lemma}

\begin{proposition}
$A(X):=
 \{A_n(X)=
{}^\vee\!X^{\otimes n}\otimes X^{\otimes n} \}_{n\in\ZZ_{\ge 0}}$
is a graded bialgebra with the following (co)multiplications:
\begin{eqnarray}
&&\mu_{m,n}:=
\mu_{{}^\vee\!X^{\otimes m}\otimes X^{\otimes m},
     {}^\vee\!X^{\otimes n}\otimes X^{\otimes n}}\,,\nonumber\\
&&\Delta _n:=\Delta_{{}^\vee\!X^{\otimes n}\otimes X^{\otimes n}}:
{}^\vee\!X^{\otimes n}\otimes X^{\otimes n}\rightarrow
({}^\vee\!X^{\otimes n}\otimes X^{\otimes n})\otimes
({}^\vee\!X^{\otimes n}\otimes X^{\otimes n})\,. \nonumber
\end{eqnarray}
Graded tensor algebra
$\T{X}$
(resp.
$\T{{}^\vee\!X}$)
is right (resp. left) $A(X)$-comodule algebra and
coalgebra.
\end{proposition}

One can carry out the same construction in the category
$\overline {\cal C}$.
The result is a bialgebra $\overline A(X)$ with
the same underlying coalgebra but with new multiplication
$\overline\mu$ where $\Psi$ is replaced by $\Psi^{-1}$.
Corresponding Hopf algebra $\overline\T{X}$ (resp. $\overline\V(X),
\overline\V^\bullet(X)$) coincide with $\T{X}_\op$ (resp. $\V(X)_\op,
\V^\bullet(X)_\op$).

\begin{figure}
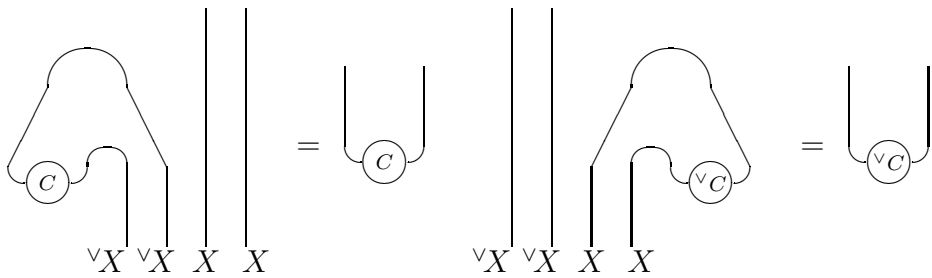

$$
\matrix{
\vvbox{\hbox{\step\coev\Step\id\step\id}
       \hbox{\dd\step\hcoev\d\step\id\step\id}
        \hbox{\coro{C}\step\id\step\id\step\id\step\id}}\cr
\Step\step\object{{}^\vee\!X}\step\object{{}^\vee\!X}\step\object{X}\step
\object{X}}
\enspace=\enspace
\vvbox{\hbox{\id\Step\id}
       \hbox{\coro{C}}}
\qquad
\matrix{
\vvbox{\hbox{\id\step\id\Step\coev}
       \hbox{\id\step\id\step\dd\hcoev\step\d}
       \hbox{\id\step\id\step\id\step\id\step\coro{{}^\vee C}}}\cr
\object{{}^\vee\!X}\step\object{{}^\vee\!X}\step\object{X}\step\object{X}\step
\Step}
\enspace=\enspace
\vvbox{\hbox{\id\Step\id}
       \hbox{\coro{{}^\vee C}}}
$$
\caption{Relations for $A(X;C)$}
\label{Figure-C-algebra}
\end{figure}

\subsec{}
Let, moreover, duality $C$ of $X$ with itself be given.
Then for each $n$ the pairing and copairing
\begin{equation}
C=C^{X^{\otimes n},X^{\otimes n}}:
X^{\otimes n}\otimes X^{\otimes n}\rightarrow \overline{1},\qquad
C=C_{X^{\otimes n},X^{\otimes n}}:
\overline{1}\rightarrow X^{\otimes n}\otimes X^{\otimes n}\,,
\label{Equation-Duality-C} \end{equation}
described by the diagram on Fig.\ref{Fig-Pairing}b and
by the input-output reversed diagram, define duality of $X^{\otimes n}$
with itself.
Let as define pairing
${}^{\vee}C=:{}^\vee\!X^{\otimes n}\otimes
{}^\vee\!X^{\otimes n}\rightarrow \overline{1} $
as morphism left dual to copairing $C_{X^{\otimes n},X^{\otimes n}}$.
We denote by $A(X;C)$ filtered algebra which is a factor-algebra of $A(X)$ by
the
ideal 'generated by relations' on Fig.\ref{Figure-C-algebra} which
means that pairings $C$ and ${}^\vee C$ are invariant with respect to
coactions of $A(X)$ on $\T{X}$ and $\T{{}^\vee\!X}$ respectively.
And let $\T{X;C}$ be factor-algebra of $\T{X}$ by relations
$C^{{}^\vee\!X\otimes X}=\underline 1$.

\begin{proposition}
$A(X;C)$ is a braided group with antipode and its inverse
given by the diagram on Fig\ref{Figure-CC-algebra}a.
$\T{X;C}$ is a right comodule algebra over
$A(X;C)$.
\end{proposition}
\begin{figure}
$$
\matrix{
\matrix{
\object{{}^\vee\!X^{\otimes n}}\Step\Step\object{X^{\otimes n}}\cr
\vvbox{\hbox{\id\step\ro{C}\step\id}
       \hbox{\hev\hcoev\step\hx}
       \hbox{\dd\step\hx\step\id}
       \hbox{\coro{C}\step\id\step\id}}}
\qquad\quad
\matrix{\object{{}^\vee X^{\otimes n}}\Step\object{X^{\otimes n}}
        \Step\Step\Step\cr
\vvbox{\hbox{\id\step\id\Step\hcoev\step\ro{C}}
       \hbox{\d\coro{C}\step\id\step\x}
       \hbox{\step\d\Step\hx\Step\id}
       \hbox{\Step\ev\step\id\Step\id}}\cr
\step\Step\Step\object{{}^\vee X^{\otimes n}}\Step\object{X^{\otimes n}}}
&\enspace&
\matrix{\object{{}^\vee\!X^{\otimes m}}\Step\object{X^{\otimes m}}
        \Step\object{{}^\vee\!X^{\otimes n}}\Step\object{X^{\otimes
n}}\cr \vvbox{\hbox{\id\step\hx\step\id}
       \hbox{\id\step\id\step\hx\obj{R}}
       \hbox{\d\hev\dd}
       \hbox{\step\hev}}}
\cr
\begin{picture}(.5,.5)\end{picture}&&\cr
\hbox{\scriptsize a) antipode and its inverse  for
$A(X;C)$}&&
\hbox{\scriptsize b)
Coquasitriangular structure on $A(X,R)$ } }
$$
\caption{}
\label{Figure-CC-algebra}
\end{figure}

Let $I=\{ I_n\in{}^\vee\!X^{\otimes n}\otimes X^{\otimes n}\}$
be a graded ideal of algebra $A(X)$ generated by relations
\begin{equation}
\label{Eq-FRT-rel}
{}^\vee\!R\otimes\id_{ X^{\otimes 2}}-
\id_{{}^\vee\!X^{\otimes 2}}\otimes R:\,
{}^\vee\!X^{\otimes 2}\otimes X^{\otimes 2}\rightarrow
{}^\vee\!X^{\otimes 2}\otimes X^{\otimes 2}
X^{\otimes 2}\otimes X^{\otimes 2}
\end{equation}
or explicitly
\begin{equation}
I_n=\bigcup_{i=1}^{n-1}({}^\vee\!X^{\otimes (i-1)}\otimes {}^\vee\!R\otimes
{}^\vee\!X^{\otimes (n-i-1)}\otimes X^{\otimes n}\rightarrow
{}^\vee\!X^{\otimes n}\otimes X^{\otimes (n-i-1)}\otimes R\otimes
X^{\otimes (i-1)} )
\end{equation}

\begin{lemma}
\label{Lemma-AXR}
Ideal $I$ described above is a biideal of $A(X)$.
We denote by $A(X,R)$ corresponding factor-bialgebra.
\end{lemma}

Let, moreover, $C$ be morphism in $\C(X,{}^\vee\!X;R)$,
i.e. pairing $C$ 'flows' along braids labeled by $R$.
Then we define a bialgebra $A(X,R;C)$ which is a factor-algebra of $A(X)$
by an ideal generated by both (\ref{Eq-FRT-rel}) and the relations given on
Fig.{\ref{Figure-C-algebra}.
Similarly, one can define the factor-algebras
$\overline{A}(X,R^{-1})$ and $\overline{A}(X,R^{-1};C)$ the
bialgebra $\overline{A}(X)$ in the category
$\overline{\cal C}$.

\begin{lemma}
A family of pairings
$\rho_{m,n}:=({}^\vee\!X^{\otimes m}\otimes X^{\otimes m})\otimes
({}^\vee\!X^{\otimes m}\otimes X^{\otimes m})\rightarrow \underline{1}$
described by the diagram on Fig.\ref{Figure-CC-algebra}b define
bialgebra pairings
\begin{eqnarray}
&&\rho_{A(X,R)}:\overline{A}(X,R^{-1})^{\rm op}\otimes A(X,R)\rightarrow
\overline{1}\,,\nonumber\\
&&\rho:_{A(X,R;C)}:\overline{A}(X,R^{-1};C)^{\rm op}\otimes A(X,R;C)
\rightarrow
\overline{1}\,.
\end{eqnarray}
\end{lemma}

\begin{theorem}
$(A(X,R),\overline{A}(X,R^{-1}),\rho)$ is a coquasitriangular
bialgebra and its factor-algebra $(A(X,R;C),\overline{A}(X,R^{-1};C),\rho)$
is a dual quantum braided group in $\C$ (or, more precisely, in a certain
category of 'graded spaces' over $\C$).
'Second inverse' $\rho^\sim=\{\rho^\sim_{m,n}\}$ to quasitriangular
structure $\rho$ takes the form:
\begin{displaymath}
\rho^\sim_{m,n}:=(\cup\otimes\cup)\circ
({}^\vee\!X^{\otimes  m}\otimes
 \Psi^{-1}_{X^{\otimes m},{}^\vee\!X^{\otimes n}}\otimes X^{\otimes n} )
\end{displaymath}
\end{theorem}

\begin{lemma}
$X^{\otimes n}$ equipped with coaction (\ref{Delta-r}) is an
object of ${\cal C}^{\cO{A(X,R),\overline A(X,R^{-1})}}$.
Braiding $\Psi_{X^{\otimes m},X^{\otimes n}}$ in this category equals to
$R_{X^{\otimes m},X^{\otimes n}}$.
Corresponding right action (defined by the first diagram on
Fig.\ref{Fig-Jurco}b)
\begin{displaymath}
\mu^{X^{\otimes n}}_r=
\{\mu^{X^{\otimes n}}_{r,m}:
  X^{\otimes n}\otimes A(X,R)_m\rightarrow X^{\otimes n}\}
\end{displaymath}
takes the form
\begin{equation}
\mu^{X^{\otimes n}}_{r,m}=
(\cup^{{}^\vee\!X^{\otimes m},X^{\otimes m}}\otimes \id_{X^{\otimes n}})
\circ
(\id_{{}^\vee\!X^{\otimes m}}\otimes R_{X^{\otimes n},X^{\otimes m}})\circ
(\Psi_{X^{\otimes n},{}^\vee\!X^{\otimes m}}\otimes\id_{X^{\otimes m}})\,.
\end{equation}
\end{lemma}

One can carry out constructions from \ref{tensor-Hopf}
for $X\in\Obj({\cal C}^{\cO{A(X,R),\overline A(X,R^{-1})}})$
to get the Hopf algebras $\T{X,R},\V(X,R),\V^\bullet(X,R)$ in
${\cal C}^{\cO{A(X,R),\overline A(X,R^{-1})}}$,
where letter '$R$' is added to specify a category.
A pair $(\V(X,R),\V(X,R)^\op)$ is a quantum braided group in
${\cal C}^{\cO{A(X,R),\overline A(X,R^{-1})}}$ with the trivial
coquasitriangular structure
$\epsilon^{\V(X,R)}\otimes\epsilon_{\V(X,R)}$.
Generalized bosonization theorem \cite{Bespalov2} allows us to define
a quantum braided group
$(A(X,R)\ltimes\V(X,R),
  \overline A(X,R^{-1})\ltimes\V(X,R)^\op)$
with coquasitriangular structure
$(\id_{A(X,R)}\otimes\epsilon_{\V(X,R)}\otimes
           \id_{A(X,R)}\otimes\epsilon_{\V(X,R)})\circ\rho_{A(X,R)}$
in $\C$
which is an analog of algebra of functions on inhomogeneous linear group.
The same construction performed for algebras $\overline A(X,R^{-1})$ and
$\overline\V(X,R^{-1})=\V(X,R)_\op$ produces quantum braided group in
$\overline\C$.
But in this way we obtain another corresponding quantum braided group
in $\C$.

\begin{figure}
$$
  \matrix{
  \matrix{ \vvbox{\hbox{\coev}
       \hhbox{\krl\id\step\ld}}\cr
\object{X}\step\object{A}\step\object{{}^\vee X}}
\enspace=\enspace
\matrix{
\vvbox{\hbox{\coev}
       \hhbox{\krl\rd\step\id}}\cr
\object{X}\step\object{A}\step\object{{}^\vee X}}
&\qquad&
\matrix{\object{X}\hstep\Step\object{A}\cr
\vvbox{\hhbox{\krl\hrd\Step\id}
       \hbox{\id\hstep\coro{\rho}}}}
\quad
\matrix{\object{A}\Step\object{{}^\vee X}\cr
\vvbox{\hbox{\step\x}
       \hbox{\ld\Step\id}
       \hbox{\hxx\Step\id}
       \hbox{\id\step\coro{\rho}}}}
&\qquad&
\matrix{
\vvbox{\hbox{\coev}
       \hhbox{\krl\rd\step\id}
       \hbox{\hxx\step\id}
       \hbox{\O{S^-}\step\hxx}
       \hbox{\hxx\step\id}}\cr
\object{{}^\vee X}\step\object{A}\step\object{X}}
\cr
\hbox{\scriptsize a) coaction $\Delta^{{}^\vee X}_\ell$ }
&&
\hbox{\scriptsize b) actions $\mu^X_r,\;\mu^{{}^\vee X}_\ell$ }
&&
\hbox{\scriptsize c) bi-invariant $\omega$ }
}
$$
\caption{}
\label{Fig-Jurco}
\end{figure}

\subsec{}
Let $(A,\overline A,\rho)$ be a dual quantum braided group in $\cal C$,
$(X,\Delta^X_r)\in{\rm Obj}({\cal C}^{\cO{A,\overline A}})\,,$
${}^\vee X$ left dual to $X$ in $\cal C$ with left comodule structure
$\Delta^{{}^\vee X}_\ell$ defined by the condition in
Fig.\ref{Fig-Jurco}a.
Then $X$ (resp. ${}^\vee X$) equipped with right (resp. left)
$A$-module structure as shown in Fig.\ref{Fig-Jurco}b becomes right
(resp. left) crossed module over $A$.
According to general theory \cite{BD1} the object
$\Gamma:={}^\vee X\otimes A\otimes X$ with actions and coactions
\begin{eqnarray}
&&\mu^\Gamma_\ell:=
(\mu^{{}^\vee X}_\ell\otimes m_A)\circ
(A\otimes\Psi_{A,{}^\vee X}\otimes A)\circ
(\Delta_A\otimes{}^\vee X\otimes A)
\otimes X\,,\nonumber\\
&&\mu^\Gamma_r:=
{}^\vee X\otimes
(m_A\otimes\mu^X_r)\circ
(A\otimes\Psi_{X,A}\otimes A)\circ
(A\otimes X\otimes\Delta_A)\,,
\nonumber\\
&&\Delta^\Gamma_\ell:=
(m_A\otimes{}^\vee X\otimes A)\circ
(A\otimes\Psi_{{}^\vee X,A}\otimes A)\circ
(\Delta^{{}^\vee X}_\ell\otimes\Delta_A)
\otimes X\,,
\nonumber\\
&&\Delta^\Gamma_r:=
{}^\vee X\otimes
(A\otimes X\otimes m_A)\circ
(A\otimes\Psi_{A,X}\otimes A)\circ
(\Delta_A\otimes\Delta^X_\ell)\,.
\nonumber
\end{eqnarray}
is a Hopf bimodule over $A$.
Morphism $\omega$ defined on Fig.\ref{Fig-Jurco}c is a bicomodule
morphism where $\underline 1$ is equipped with trivial left and right
actions equal to $\epsilon_A$.
'Commutant with $\omega$':
\begin{equation}
{\rm d}:=
\mu^\Gamma_r\circ (\omega\otimes A)-
\mu^\Gamma_\ell\circ(A\otimes\omega):
A\rightarrow\Gamma
\end{equation}
is a first order bicovariant derivative in sense of Woronowicz
\cite{Wor} (See \cite{BD} for a fully braided context).

In our case $A=A(X,R;C)$ 'biinvariant' $\omega$ equals to
$\Psi^{-1}_{{}^\vee\!X,{}^\vee\!X}\circ{}^\vee\!C\otimes C:
  \underline 1\rightarrow{}^\vee\!X^{\otimes 2}\otimes X^{\otimes 2}$.

\end{document}